\pdfoutput=1
\documentclass[12pt,leqno]{amsart}
\usepackage[letterpaper,margin=1in]{geometry}
\usepackage[T1]{fontenc}
\usepackage{lmodern}
\usepackage{amsmath,amssymb,mathtools}
\usepackage{microtype}
\DeclareMicrotypeSet{papertext}{encoding=T1,family={rm*,sf*}}
\UseMicrotypeSet[expansion]{papertext}
\usepackage[unicode,bookmarksnumbered]{hyperref}
\numberwithin{equation}{section}
\newtheorem{theorem}{Theorem}[section]
\newtheorem{proposition}[theorem]{Proposition}

\newtheorem{lemma}[theorem]{Lemma}
\theoremstyle{plain}

\theoremstyle{remark}
\newtheorem{remark}[theorem]{Remark}
\newcommand{\PP}{\mathbb P}
\newcommand{\EE}{\mathbb E}
\newcommand{\RR}{\mathbb R}
\newcommand{\NN}{\mathbb N}
\newcommand{\HH}{\mathcal H}
\newcommand{\FF}{\mathcal F}
\newcommand{\GG}{\mathcal G}
\newcommand{\one}{\mathbf 1}
\newcommand{\Ball}{\mathbb B}
\newcommand{\dd}{\,\mathrm d}
\DeclareMathOperator{\diam}{diam}
\DeclareMathOperator{\dist}{dist}
\hypersetup{
  colorlinks=false,
  pdfborder={0 0 1},
  citebordercolor={0 1 0},
  pdftitle={Hausdorff measures of planar Brownian level sets},
  pdfsubject={Uniform gauge classification and the critical mass spectrum},
  pdfauthor={Yingyu Zhao, Yunhao Zhao}
}

\title{Hausdorff measures of planar Brownian level sets}
\date{}
\author{Yingyu Zhao}
\address{Yingyu Zhao, Qiuzhen College, Tsinghua University, 1 Qinghuayuan, Haidian District, 100084 Beijing, China}
\email{zhao-yy25@mails.tsinghua.edu.cn}
\author{Yunhao Zhao}
\address{Yunhao Zhao, School of Advanced Interdisciplinary Sciences, University of Chinese Academy of Sciences, 1 Yanqihu East Road, Huairou District, 101408 Beijing, China}
\email{zhaoyunhao23@mails.ucas.ac.cn}

\begin{document}

\begin{abstract}
For any Hausdorff gauge \(h\), we prove that all level sets of planar Brownian motion have zero \(h\)-measure almost surely if and only if \(\liminf_{r\to0^+}h(r)\log(1/r)=0\). This resolves Taylor's question for arbitrary Hausdorff gauges. The proof yields \(0<\HH^{h_0}(L_x)<\infty\) for some \(x\) almost surely, where \(h_0(r)=1/\log(1/r)\). The Hausdorff dimensions of the spatial mass superlevel sets form a deterministic, strictly decreasing bi-Lipschitz spectrum. For each fixed positive mass below the supremum, the exact-mass set almost surely has the same dimension as the corresponding superlevel set.
\end{abstract}

\maketitle
\flushbottom

\section{Introduction}\label{sec:intro}

Planar Brownian motion almost surely visits some points uncountably many times \cite{DEK,LeGall}. Kaufman's dimension-doubling theorem \cite{Kaufman} shows that every level set has Hausdorff dimension zero on a single event of probability one; see also \cite[Corollary~9.29]{MP}. The size of these exceptional time sets is therefore a question about Hausdorff gauges.

Let \(B=(B_t)_{t\geqslant0}\) be planar Brownian motion started at the origin with generator \(\Delta/2\), where \(\Delta\) is the Laplacian on \(\RR^2\). We write
\[
 L_x(I)=\{t\in I:B_t=x\},\quad L_x=L_x([0,\infty)).
\]
Taylor \cite[Problem~5]{Taylor} asked for a characterization of the gauges \(h\) for which
\[
 \PP\{\HH^h(L_x)=0\text{ for every }x\in\RR^2\}=1.
\]
For \(h(r)=(\log(1/r))^{-\alpha}\), the results of Perkins and Taylor \cite{PT87} imply uniform vanishing when \(\alpha>2\). They also formulated the logarithmic-power question in Problem~3.11.

Bass, Burdzy and Khoshnevisan \cite{BBK} constructed measures \(\beta_a\) for \(0<a<1/2\), each carried by points of infinite multiplicity. Using these measures, Cammarota and M\"orters \cite[Theorem~1.1]{CM} proved for each fixed \(a\in(0,1/2)\) that \(\HH^{h_*}(L_x)>0\) almost surely for \(\beta_a\)-almost every \(x\), where
\[
 h_*(r)=\frac{\log\log\log(1/r)}{\log(1/r)}.
\]
For each gauge \(h\) satisfying \(h(r)\log(1/r)\to0\) as \(r\to0^+\), they proved \(\HH^h(L_x)=0\) simultaneously for all \(x\) with probability one \cite[Theorem~1.2]{CM}. These bounds settle the logarithmic-power question for \(0<\alpha<1\) and \(\alpha>1\), respectively. For general gauges, they conjectured uniform vanishing whenever \(h=o(h_*)\) \cite[Section~4]{CM}. Huang \cite{Huang} obtained related results for the logarithmic Minkowski dimension of visit-time sets.

Our first theorem gives the classification for arbitrary gauges. The main step is to construct level sets of positive measure at the gauge
\begin{equation}\label{eq:h0}
 h_0(0)=0,\quad h_0(r)=\frac1{\log(1/r)},\quad 0<r<1.
\end{equation}
A gauge is a nondecreasing, right-continuous function \(h\) on \([0,r_h)\) with \(h(0)=0\), \(h(r)>0\) for \(r>0\), and \(h(r)\to0\) as \(r\to0^+\). Denote this class by \(\GG\). The Hausdorff \(h\)-content at scale \(0<\delta<r_h\) and the corresponding measure are
\begin{equation}\label{eq:hausdorff}
 \begin{aligned}
 \HH^h_\delta(E)&=\inf\left\{\sum_i h(\diam U_i):
 E\subseteq\bigcup_iU_i,\ \diam U_i<\delta\right\},\\
 \HH^h(E)&=\lim_{\delta\to0^+}\HH^h_\delta(E),
 \end{aligned}
\end{equation}
where the infimum is over countable covers of \(E\). The class \(\mathcal I\) consists of all nondegenerate intervals in \([0,\infty)\), including unbounded intervals. Write
\begin{equation}\label{eq:Lambda}
 \Lambda(h)=\liminf_{r\to0^+}h(r)\log(1/r),\quad h\in\GG.
\end{equation}

\begin{theorem}\label{thm:classification}
There exist deterministic constants \(0<m_B\leqslant C_0<\infty\) and an event \(\Omega_*\) of probability one. On \(\Omega_*\), the following hold for all \(h\in\GG\) and \(I\in\mathcal I\).
\begin{enumerate}
\renewcommand{\labelenumi}{(\roman{enumi})}
\item \(\HH^h(L_x(I))=0\) for every \(x\in\RR^2\) if and only if \(\Lambda(h)=0\).
\item \(\displaystyle\sup_x\HH^{h_0}(L_x(I))=m_B\).
\item \(\displaystyle m_B\Lambda(h)\leqslant\sup_x\HH^h(L_x(I)) \leqslant C_0\Lambda(h),\quad\Lambda(h)<\infty.\)
\item If \(\Lambda(h)=\infty\), then \(\HH^h(L_x(I))=\infty\) for some \(x\in\RR^2\).
\end{enumerate}
\end{theorem}

\begin{remark}
Since \(h_0=o(h_*)\), part~(ii) disproves the conjecture in \cite[Section~4]{CM}.
\end{remark}

We next study the spatial distribution of \(f_I(x)=\HH^{h_0}(L_x(I))\). For \(m\geqslant0\), write
\[
 \mathcal S_m(I)=\{x:f_I(x)>m\},\quad
 \mathcal T_m(I)=\{x:f_I(x)\geqslant m\}.
\]
We use the convention \(\dim_{\mathrm H}\varnothing=0\).

\begin{theorem}\label{thm:mass-spectrum}
There exist absolute constants \(c_0,C_1>0\) and a deterministic function \(D:[0,m_B]\to[0,2]\) with \(D(0)=2\) and \(D(m_B)=0\) such that
\begin{equation}\label{eq:spectrum-lipschitz}
 c_0(m_2-m_1)\leqslant D(m_1)-D(m_2)\leqslant C_1(m_2-m_1),
 \quad 0\leqslant m_1<m_2\leqslant m_B.
\end{equation}
The event \(\Omega_*\) can be chosen so that every \(I\in\mathcal I\) satisfies
\begin{equation}\label{eq:spectrum-deterministic}
 \mathcal S_m(I)\ne\varnothing,\quad
 \dim_{\mathrm H}\mathcal S_m(I)=D(m),\quad 0\leqslant m<m_B,
\end{equation}
and \(\dim_{\mathrm H}\mathcal T_m(I)=D(m)\) for \(0<m\leqslant m_B\). For every open set \(U\) meeting \(B(I)\),
\[
 \dim_{\mathrm H}(\mathcal S_m(I)\cap U)=D(m),
 \quad 0\leqslant m<m_B.
\]
\end{theorem}

\begin{remark}
Dembo, Peres, Rosen and Zeitouni \cite[Theorem~1.3]{DPRZ} obtained the spectrum \(2-a\), \(0\leqslant a\leqslant2\), for occupation thickness on a fixed finite time interval, with normalization \(r^2\log^2(1/r)\). A\"id\'ekon, Hu and Shi \cite{AHS} extended the measures of \cite{BBK} to all subcritical thickness parameters. For each fixed parameter, their Brownian excursion construction almost surely gives continuous local times at almost every center with respect to the corresponding measure; the associated Stieltjes measures have the visit-time sets as their support \cite[Theorem~1.2]{AHS}. Related Brownian multiplicative chaos measures were constructed and studied by Jego \cite{Jego,JegoCritical,JegoCharacterisation}. Our parameter is the Hausdorff measure of the time set \(L_x(I)\); our decomposition uses the finite, nonatomic measure \(\HH^{h_0}|_{L_x}\) and its additivity in time.
\end{remark}

For \(m>0\), write
\[
 \mathcal X_m(I)=\{x:f_I(x)=m\}.
\]
\begin{theorem}\label{thm:exact-mass-intro}
For each fixed \(m\in(0,m_B)\), there is an event of probability one on which
\[
 \mathcal X_m(I)\ne\varnothing,\quad
 \dim_{\mathrm H}\mathcal X_m(I)=D(m)
\]
for every \(I\in\mathcal I\). Moreover, for every such \(I\) and every open set \(U\) meeting \(B(I)\), \(\dim_{\mathrm H}(\mathcal X_m(I)\cap U)=D(m)\).
\end{theorem}

\begin{remark}
Theorem~\ref{thm:exact-mass-intro} holds on one event for any prescribed countable family of masses. A common event for all \(m\in(0,m_B)\) has not been established. The exact value of \(m_B\), the form of \(D\), and whether \(m_B\) is attained remain open.
\end{remark}

Delayed-return covers as in \cite{CM} give the upper bound on one event for all gauges. For the lower bound, we adapt the nested return construction of \cite{DEK,MP}, using intervals of deterministic lengths \(r_j=e^{-a2^j}\) with quantitative separation. Equal division of mass assigns weight \(2^{-j}=a h_0(r_j)\) at generation \(j\). By separation, the resulting measure \(\mu\) satisfies \(\mu(J)\leqslant2a h_0(|J|)\) for every sufficiently small interval \(J\). This controls covers with varying diameters. A two-point estimate yields finite-energy measures on the successful centers. Brownian scaling and independence on half-intervals make the mass supremum deterministic.

For the spectrum, we split the mass between past and future. Conditional on the past, upper and lower estimates for the future path on a fixed spatial target give \eqref{eq:spectrum-lipschitz} in Section~\ref{sec:mass-spectrum}. In Section~\ref{sec:exact-mass}, an optional section selects a stopping time at which the accumulated mass equals \(m\). Point polarity preserves that mass. We then use independent percolation sets to transfer the hitting statement to Hausdorff dimension, following the method of \cite[Lemma~3.4 and Remark~3.5]{KPX}. This argument applies to each fixed \(h\in\GG\) with \(\Lambda(h)<\infty\); see Proposition~\ref{thm:exact-mass-transfer}.

\section{Preliminaries and uniform covers}\label{sec:gauges}

All logarithms in this paper are natural. We write \(\PP_z\) and \(\EE_z\) for the law and expectation of Brownian motion started at \(z\), and \(\Ball(x,r)\) for the open disk of radius \(r\) centered at \(x\). We equip \(\mathcal C=C([0,\infty),\RR^2)\) with the topology of uniform convergence on compact intervals. We use the completed, right-continuous natural Brownian filtration \((\FF_t)\) under the law in use.

\begin{lemma}\label{lem:gauge-properties}
Let \(h\in\GG\) and \(E\subseteq\RR\).
\begin{enumerate}
\renewcommand{\labelenumi}{(\roman{enumi})}
\item \(\displaystyle\liminf_{n\to\infty}bn\,h(e^{-bn})=\Lambda(h),\quad b>0.\)
\item \(\HH^h(E)\geqslant\Lambda(h)\HH^{h_0}(E),\quad0<\Lambda(h)<\infty.\) If \(\Lambda(h)=\infty\) and \(\HH^{h_0}(E)>0\), then \(\HH^h(E)=\infty\).
\item \(\HH^{h_0}(t+cE)=\HH^{h_0}(E),\quad c>0,\ t\in\RR.\)
\end{enumerate}
\end{lemma}
\begin{proof}
We write \(n(r)=\lceil b^{-1}\log(1/r)\rceil\). Monotonicity gives
\[
 bn(r)h(e^{-bn(r)})\leqslant h(r)\bigl(\log(1/r)+b\bigr).
\]
Since \(n(r)\to\infty\) and \(h(r)\to0\) as \(r\to0^+\),
\[
 \Lambda(h)\leqslant\liminf_{n\to\infty}bn\,h(e^{-bn})
 \leqslant\liminf_{r\to0^+}h(r)\bigl(\log(1/r)+b\bigr)
 =\Lambda(h).
\]
The first inequality follows by restricting to \(r=e^{-bn}\). This proves (i).

Fix \(0<d<\Lambda(h)\). We have \(h\geqslant dh_0\) at sufficiently small scales. Comparing cover costs and then taking \(\delta\to0^+\) gives
\[
 \HH^h_\delta(E)\geqslant d\HH^{h_0}_\delta(E),\quad
 \HH^h(E)\geqslant d\HH^{h_0}(E).
\]
Taking the supremum over \(d\) proves (ii).

Finally,
\[
 \frac{h_0(cr)}{h_0(r)}
 =\frac{\log(1/r)}{\log(1/r)-\log c}\longrightarrow1,
 \quad r\to0^+.
\]
Affine images and inverse images of covers therefore give
\[
 (1-\varepsilon)\HH^{h_0}_\delta(E)
 \leqslant\HH^{h_0}_{c\delta}(t+cE)
 \leqslant(1+\varepsilon)\HH^{h_0}_\delta(E)
\]
for every \(\varepsilon\in(0,1)\) and all sufficiently small \(\delta>0\). We take \(\delta\to0^+\) and then \(\varepsilon\to0^+\) to obtain (iii).
\end{proof}

For \(T<\infty\), \(w\in C([0,T],\RR^2)\), \(x\in\RR^2\), and \(0\leqslant a\leqslant b\leqslant T\), write
\[
 K(w,x;a,b)=\{t\in[a,b]:w(t)=x\},
\]
and denote its content and measure by
\[
 C^h_\delta(w,x;a,b)=\HH^h_\delta(K(w,x;a,b)),\quad
 f^h_{[a,b]}(w,x)=\HH^h(K(w,x;a,b)).
\]
Choose \(n_0\) so that \(1/n_0\) lies in the domain of \(h\).

\begin{lemma}\label{lem:level-set-content}
Fix \(h\in\GG\) and \(T<\infty\).
\begin{enumerate}
\renewcommand{\labelenumi}{(\roman{enumi})}
\item For every \(0<\delta<r_h\), the function \(C^h_\delta\) is jointly upper semicontinuous in \((w,x,a,b)\).
\item The functions \(f^h_{[a,b]}(w,x)\) and \(\sup_x f^h_{[a,b]}(w,x)\) are Borel in their parameters.
\item For \(m\geqslant0\),
\[
 \{x:f^h_{[a,b]}(w,x)>m\}
 =\bigcup_{n\geqslant n_0,\,k\geqslant1}
   \{x:C^h_{1/n}(w,x;a,b)\geqslant m+1/k\}.
\]
Each term is compact and has a closed graph in \((w,a,b,m,x)\).
\item The dimensions of the strict superlevels in part~(iii) and their pairwise intersections are Borel in the path, interval and threshold parameters.
\end{enumerate}
\end{lemma}
\begin{proof}
We first show that \(C^h_\delta\) can be computed using finite covers by open intervals with rational endpoints. Given \(\varepsilon>0\) and a countable cover \((U_i)\) of \(K(w,x;a,b)\) with \(\diam U_i<\delta\), right continuity of \(h\) gives open rational intervals \(J_i\supseteq U_i\) satisfying
\[
 |J_i|<\delta,\quad
 h(|J_i|)\leqslant h(\diam U_i)+\varepsilon2^{-i}.
\]
Compactness gives a finite subcover. Taking infima over the original covers and then letting \(\varepsilon\to0^+\) gives the claimed reduction.

For a fixed open \(U\subseteq\RR\), the condition \(K(w,x;a,b)\subseteq U\) is open in \((w,x,a,b)\). To see this, suppose \((w_j,x_j,a_j,b_j)\to(w,x,a,b)\) and \(t_j\in K(w_j,x_j;a_j,b_j)\setminus U\). Along a subsequence,
\[
 t_j\to t\in[a,b]\setminus U,\quad
 w(t)=\lim_j w_j(t_j)=\lim_j x_j=x.
\]
Thus the complementary condition is closed. Each strict sublevel of \(C^h_\delta\) is a union of these open conditions over finite rational covers whose total cost is below the given threshold. This proves (i).

To prove (ii), we use the identities
\[
 \begin{aligned}
 f^h_{[a,b]}(w,x)&=\sup_{n\geqslant n_0}C^h_{1/n}(w,x;a,b),\\
 \sup_x f^h_{[a,b]}(w,x)
 &=\sup_{n\geqslant n_0,\,k\geqslant1}
       \max_{|x|\leqslant k}C^h_{1/n}(w,x;a,b).
 \end{aligned}
\]
Each maximum is upper semicontinuous by (i) and compactness of \(\overline\Ball(0,k)\). Since countable suprema of Borel functions are Borel, this proves (ii). The first identity also gives the decomposition in (iii). Each term has a closed graph by (i) and is contained in the compact image \(w([a,b])\), since its threshold \(m+1/k\) is positive.

By (iii), the strict superlevel families have Borel graphs and \(\sigma\)-compact sections. The same holds for their intersections by intersecting the compact terms. Since the parameter spaces are Polish, (iv) follows from \cite[Theorem~6.1]{MM}.
\end{proof}

\subsection{A uniform bound on the whole time axis}\label{sec:upper}

Fix \(\theta\in(1/2,1)\), and choose
\[
 c_0=\log(1/\theta),\quad C_0>2/c_0.
\]
For an integer \(R\geqslant1\), write \(\tau_R=\inf\{t\geqslant0:|B_t|=R\}\). For a closed mesh square \(Q\) of side \(\varepsilon<1\), its delayed visits are
\[
 S_1=\inf\{t\geqslant0:B_t\in Q\},\quad
 S_{k+1}=\inf\{t\geqslant S_k+\varepsilon:B_t\in Q\},
\]
and \(N_Q=\#\{k:S_k<\tau_R\}\). The same notation applies to closed disks of radius \(\varepsilon\).

\begin{proposition}\label{prop:upper}
On one event of probability one, for every \(h\in\GG\) with \(\Lambda(h)<\infty\) and every \(x\in\RR^2\),
\begin{equation}\label{eq:upper}
 \HH^h(L_x)\leqslant C_0\Lambda(h).
\end{equation}
\end{proposition}
\begin{proof}
We fix \(R\) and consider the squares \(Q\) meeting \(\overline\Ball(0,R)\). Write \(z_Q\) for the center of \(Q\) and \(A=4R+4\). To estimate delayed returns to \(Q\), we use \(u_\varepsilon:\RR^2\to[0,1]\), where \(u_\varepsilon(y)\) is the probability that Brownian motion started at \(y\) hits \(\overline\Ball(z_Q,\varepsilon)\) before exiting \(\Ball(z_Q,A)\). We extend \(u_\varepsilon\) by zero outside the outer disk. The annulus hitting formula \cite[Theorem~3.18]{MP} gives
\[
 u_\varepsilon(y)=
 \begin{cases}
  1,&|y-z_Q|\leqslant\varepsilon,\\
  \log(A/|y-z_Q|)/\log(A/\varepsilon),
       &\varepsilon<|y-z_Q|<A,\\
  0,&|y-z_Q|\geqslant A.
 \end{cases}
\]
Fix \(z\in Q\cap\Ball(0,R)\). Since \(Q\subseteq\overline\Ball(z_Q,\varepsilon)\) and \(\Ball(0,R)\subseteq\Ball(z_Q,A)\), the Markov property at time \(\varepsilon\) gives
\begin{equation}\label{eq:return-potential}
 \PP_z\{B_t\in Q\text{ for some }\varepsilon\leqslant t<\tau_R\}
 \leqslant\EE_z u_\varepsilon(B_\varepsilon).
\end{equation}
For small \(\varepsilon\),
\begin{align*}
 u_\varepsilon(y)&\leqslant
 \frac{\log(A/\sqrt\varepsilon)
       +\log^+(\sqrt\varepsilon/|y-z_Q|)}{\log(A/\varepsilon)},\\
 \EE_z\log^+\frac{\sqrt\varepsilon}{|B_\varepsilon-z_Q|}
 &\leqslant\frac1{2\pi\varepsilon}
       \int_{|u|<\sqrt\varepsilon}
          \log\frac{\sqrt\varepsilon}{|u|}\dd u=1/4.
\end{align*}
Thus the supremum \(q_\varepsilon\) of these delayed-return probabilities satisfies
\begin{equation}\label{eq:return-prob}
 q_\varepsilon\leqslant
 \frac{L/2+\log A+1/4}{L+\log A}
 =1/2+O_R(L^{-1}),\quad L=\log(1/\varepsilon).
\end{equation}
For all sufficiently small \(\varepsilon\), the strong Markov property at successive \(S_k\) gives
\begin{align}
 \PP(N_Q\geqslant k)&\leqslant q_\varepsilon^{k-1}
 \leqslant\theta^{k-1},\quad k\geqslant1,
 \label{eq:geometric-tail}\\
 \PP\{N_Q\geqslant\lceil u\log(1/\varepsilon)\rceil\}
 &\leqslant\theta^{-1}\varepsilon^{c_0u},\quad u>0.
 \label{eq:count-tail}
\end{align}
These estimates also hold for closed disks of radius \(\varepsilon\).

At scale \(\varepsilon_n=2^{-n}\), the number of relevant squares is \(O_R(\varepsilon_n^{-2})\). Since \(c_0C_0>2\),
\[
 \sum_{n\geqslant n_R}
 \PP\{\exists Q:N_Q\geqslant C_0\log(1/\varepsilon_n)\}
 \leqslant C_R\sum_{n\geqslant n_R}
              \varepsilon_n^{c_0C_0-2}<\infty.
\]
Borel--Cantelli and intersection over integer \(R\) give one event of probability one on which
\[
 N_Q<C_0\log(1/\varepsilon_n)\quad\text{for every relevant }Q
\]
for every \(R\) and all sufficiently large \(n\). If \(Q\) contains \(x\), the intervals \([S_k,S_k+\varepsilon_n]\) with \(S_k<\tau_R\) cover \(L_x([0,\tau_R))\). Adding \(\{\tau_R\}\) has zero cost, so for \(2\varepsilon_n<r_h\),
\[
 \HH^h_{2\varepsilon_n}(L_x([0,\tau_R]))
 \leqslant N_Qh(\varepsilon_n)
 \leqslant C_0h(\varepsilon_n)\log(1/\varepsilon_n).
\]
Taking lower limits and applying Lemma~\ref{lem:gauge-properties}(i) gives
\[
 \HH^h(L_x([0,\tau_R]))
 \leqslant C_0\liminf_n h(\varepsilon_n)\log(1/\varepsilon_n)
 =C_0\Lambda(h)
\]
simultaneously for all \(h,x,R\). Continuity from below as \(R\to\infty\) gives \eqref{eq:upper}.
\end{proof}

\section{Critical levels from return trees}\label{sec:node}

Binary nesting of visit times appears in the work of Dvoretzky, Erd\H{o}s and Kakutani \cite[p.~179, step~(G)]{DEK}. M\"orters and Peres \cite[proof of Theorem~9.24]{MP} construct uncountable level sets from nested time intervals. Their argument uses one- and two-point estimates for Brownian excursions \cite[Lemma~9.27]{MP}. We use binary trees of nested time intervals with deterministic lengths and quantitative separation. For \(j\geqslant0\),
\begin{equation}\label{eq:rj}
 r_j=e^{-a2^j},\quad r_{j+1}=r_j^2,\quad
 h_0(r_j)=a^{-1}2^{-j},\quad a>1.
\end{equation}
Equal division of mass and separation give an \(h_0\)-Frostman bound.

For a finite Borel measure \(\nu\) on \(\RR^2\), write
\[
 \mathcal E_s(\nu)=\iint|x-y|^{-s}\nu(\dd x)\nu(\dd y),\quad s>0.
\]

\subsection{Return events}

We define the return events on the path space \(\mathcal C\), with coordinate process \(B_t(w)=w(t)\). For \(w\in\mathcal C\), write \((\theta_t w)(u)=w(t+u)\) and
\[
 H_\rho(w,x)=\inf\{t\geqslant0:|w(t)-x|=\rho\},\quad\rho>0.
\]
For a stopping time \(S\), we write \(\theta_SB=(B_{S+u})_{u\geqslant0}\) on \(\{S<\infty\}\). We take \(\inf\varnothing=\infty\) and declare every node with an infinite starting time unsuccessful. Fix \(M=4\), and define
\begin{equation}\label{eq:G0}
 G_0(x,r)=\left\{w\in\mathcal C:
 \sup_{0\leqslant t\leqslant r}|w(t)-x|\leqslant M\sqrt r,\quad
 |w(r)-x|\leqslant\sqrt r\right\}.
\end{equation}
The endpoint bound enters the recursive return estimate; spatial confinement is used in the two-point bound. For \(j\geqslant1\), define recursively
\begin{gather}
 \begin{aligned}
 \sigma_1(x,r)&=\inf\{t\geqslant0:|w(t)-x|=r^2\},\\
 \sigma_2(x,r)&=\inf\{t\geqslant\sigma_1(x,r)+2r^2:|w(t)-x|=r^2\},
 \end{aligned}\label{eq:sigmas}\\
 G_j(x,r)=G_0(x,r)\cap\{\sigma_1\leqslant r/8,\ \sigma_2\leqslant r/2\}
 \cap\bigcap_{i=1}^2\{\theta_{\sigma_i}B\in G_{j-1}(x,r^2)\}.
 \label{eq:Gj}
\end{gather}
On \(G_j\), the intervals \([\sigma_i,\sigma_i+r^2]\) lie in \([0,r]\) and are separated by at least \(r^2\) whenever \(r\leqslant1/2\).

The identity
\[
 \{H_\rho\leqslant t\}
 =\left\{\min_{0\leqslant u\leqslant t}
            \bigl||w(u)-x|-\rho\bigr|=0\right\}
\]
and continuity of the shift at finite times give joint Borel measurability of all starting times and events in \((w,x)\). For \(r\leqslant1/2\), the bound \(r/2+r^2\leqslant r\) and induction give
\[
 G_j(x,r)\in\FF_r,\quad G_{j+1}(x,r)\subseteq G_j(x,r).
\]
By translation and rotation invariance, the probability
\begin{equation}\label{eq:pj}
 p_j(r)=\PP_z(G_j(x,r)),\quad |z-x|=r,
\end{equation}
depends only on \(j\) and \(r\).

\begin{lemma}\label{lem:one-node}
There are \(r_*>0\) and \(b,q\in(0,1)\) such that
\begin{equation}\label{eq:one-node}
 p_0(r)\geqslant b,\quad p_j(r)\geqslant q\,p_{j-1}(r^2)^2
\end{equation}
for \(0<r\leqslant r_*\) and \(j\geqslant1\).
\end{lemma}
\begin{proof}
We write \(R_r=\sqrt r/16\). Let \(T_u\) be the first hitting time of the circle of radius \(u\) about \(x\). When \(2r<R_r\), the annulus formula and the mean exit time give
\[
 \begin{aligned}
 \PP_z(T_{r^2}<T_{R_r})
 &\geqslant\frac{\log(R_r/(2r))}{\log(R_r/r^2)}\longrightarrow1/3,
 \quad r\to0^+,\\
 \PP_z(T_{R_r}>r/8)
 &\leqslant\frac8r\EE_zT_{R_r}
 =\frac{4(R_r^2-|z-x|^2)}r\leqslant1/64,
 \end{aligned}
\]
uniformly for \(r^2\leqslant|z-x|\leqslant2r\). The exit-time estimate also holds for \(|z-x|<r^2\), where continuity gives \(\PP_z(T_{r^2}<T_{R_r})=1\). We choose \(r_*>0\) so that \(2r<R_r\) and the logarithmic ratio above is at least \(1/4\) whenever \(0<r\leqslant r_*\). Excluding the event \(\{T_{R_r}>r/8\}\) then gives
\begin{equation}\label{eq:short-hit}
 \begin{aligned}
 \PP_z\{T_{r^2}\leqslant r/8,\ T_{r^2}<T_{R_r}\}
 &\geqslant\PP_z(T_{r^2}<T_{R_r})-\PP_z(T_{R_r}>r/8)\\
 &\geqslant1/4-1/64\geqslant1/8
 \end{aligned}
\end{equation}
uniformly for \(|z-x|\leqslant2r\). We use the confinement probabilities
\[
 g=\PP_0\{\sup_{0\leqslant t\leqslant1}|B_t|\leqslant1\},\quad
 b=\PP_0\{\sup_{0\leqslant t\leqslant1}|B_t|\leqslant1/2\}.
\]
Both are positive by \cite[Exercise~1.8]{MP} and independence of the coordinates. With \(|z-x|=r\leqslant\sqrt r/2\), the triangle inequality and Brownian scaling give
\[
 \begin{aligned}
 p_0(r)&\geqslant\PP_z\{\sup_{0\leqslant t\leqslant r}|B_t-x|\leqslant\sqrt r\}\\
 &\geqslant\PP_z\{\sup_{0\leqslant t\leqslant r}|B_t-z|\leqslant\sqrt r/2\}=b.
 \end{aligned}
\]

To prove the recursive estimate for \(j\geqslant1\), we require the event in \eqref{eq:short-hit}, followed by \(G_{j-1}(x,r^2)\) at \(S=\sigma_1\). At \(V=S+r^2\), \(|B_V-x|\leqslant r\). Write \(W=S+2r^2\), so \(W-V=r^2\). The strong Markov property at \(V\) and Brownian scaling give
\[
 \PP_z\!\left(\sup_{V\leqslant t\leqslant W}|B_t-B_V|\leqslant r
        \mid\FF_V\right)
 =\PP_0\{\sup_{0\leqslant t\leqslant r^2}|B_t|\leqslant r\}=g.
\]
On this event, \(|B_W-x|\leqslant2r\), so another application of \eqref{eq:short-hit} gives a return at \(T=\sigma_2\) satisfying
\[
 T\leqslant W+r/8\leqslant r/4+2r^2\leqslant r/2.
\]
We then require \(G_{j-1}(x,r^2)\) at \(T\). From \(U=T+r^2\), we also require
\[
 \sup_{U\leqslant t\leqslant r}|B_t-B_U|\leqslant\sqrt r/2.
\]
On these events, \(r-U\in[0,r]\) is \(\FF_U\)-measurable, so the last event has conditional probability at least \(b\). For sufficiently small \(r_*\), these requirements confine the path within distance
\[
 \max\{R_r,Mr,2r,r+\sqrt r/2\}\leqslant M\sqrt r
\]
of \(x\), and its final position within \(r+\sqrt r/2\leqslant\sqrt r\). They therefore imply \(G_j\). Successive conditioning gives
\[
 p_j(r)\geqslant gb\,p_{j-1}(r^2)^2/64
 =q\,p_{j-1}(r^2)^2,\quad q=gb/64.\qedhere
\]
\end{proof}

We record two consequences of the recursion for the two-point estimate in Lemma~\ref{lem:fixed-forest}(iii). For \(e^{-a}\leqslant r_*\), iteration of \eqref{eq:one-node} gives
\begin{align}
 p_n(r_0)&\geqslant q^{2^k-1}p_{n-k}(r_k)^{2^k},
 \quad 0\leqslant k\leqslant n,\label{eq:iterate}\\
 p_n(r_0)&\geqslant q^{2^n-1}b^{2^n}>0.\label{eq:crude}
\end{align}

\subsection{Products and spatial estimates}

\begin{lemma}\label{lem:random-order-windows}
Let \(S_1,\ldots,S_m\) be almost surely finite stopping times, \(\ell>0\), and \(A_i\in\FF_{S_i+\ell}\). If \(\EE(\one_{A_i}\mid\FF_{S_i})\leqslant p\) for every \(i\), then
\[
 \PP\left(\bigcap_{i=1}^m A_i,\
     [S_i,S_i+\ell]\text{ have disjoint interiors}\right)\leqslant p^m.
\]
\end{lemma}
\begin{proof}
We order \(S_1,\ldots,S_m\) increasingly, breaking ties by the original indices. Let \(J_k\) be the original index of the \(k\)th time in this order, and write \(T_k=S_{J_k}\). Then \(T_k\) is a stopping time, since
\[
 \{T_k\leqslant t\}=\{\#\{i:S_i\leqslant t\}\geqslant k\}.
\]
The rank of \((S_i,i)\) is determined at \(S_i\), so \(\{J_k=i\}\in\FF_{S_i}\cap\FF_{T_k}\). Hence \(D_k=A_{J_k}\in\FF_{T_k+\ell}\) and
\[
 \EE(\one_{D_k}\mid\FF_{T_k})
 =\sum_i\one_{\{J_k=i\}}\EE(\one_{A_i}\mid\FF_{S_i})\leqslant p.
\]
For \(E_1=D_1\) and \(E_k=E_{k-1}\cap\{T_{k-1}+\ell\leqslant T_k\}\cap D_k\), induction gives \(E_k\in\FF_{T_k+\ell}\), so
\[
 E_{k-1}\cap\{T_{k-1}+\ell\leqslant T_k\}\in\FF_{T_k}.
\]
Therefore
\[
 \PP(E_k)
 =\EE\!\left[\one_{E_{k-1}\cap\{T_{k-1}+\ell\leqslant T_k\}}
              \EE(\one_{D_k}\mid\FF_{T_k})\right]
 \leqslant p\PP(E_{k-1})\leqslant p^k,
\]
and \(E_m\) is the event in the statement.
\end{proof}

We now combine these trees in separated time intervals so that the \(h_0\)-measures of the limiting sets add. Choose absolute constants
\[
 K\geqslant\max\{q^{-2},(qb)^{-1},2\},\quad K_0=4\log K,
 \quad e^{-a_0}\leqslant r_*,\quad e^{-a_0}<1/4.
\]
Fix \(a\geqslant a_0\) and an integer \(k\geqslant1\), and write \(\eta=K_0k/a\) for the exponent in the two-point estimate. We write \(|v|\) for the length of a finite binary word \(v\), with \(|\varnothing|=0\). The \(k\) trees start at
\[
 R_1(x)=H_{r_0}(B,x),\quad
 R_{i+1}(x)=\inf\{t\geqslant R_i(x)+r_0+1:|B_t-x|=r_0\}.
\]
Within the \(i\)th tree, we write \(S_{i,\varnothing}(x)=R_i(x)\) and define recursively for \(v\in\{0,1\}^j\)
\[
 \begin{aligned}
 S_{i,v0}(x)&=\inf\{t\geqslant S_{i,v}(x):|B_t-x|=r_{j+1}\},\\
 S_{i,v1}(x)&=\inf\{t\geqslant S_{i,v0}(x)+2r_{j+1}:|B_t-x|=r_{j+1}\}.
 \end{aligned}
\]
For each fixed \(x\), these are almost surely finite stopping times. The event of success through generation \(n\) in all \(k\) trees with bounded waiting times is
\[
 F_{n,k}(x)=\{R_1\leqslant1\}
 \cap\bigcap_{i<k}\{R_{i+1}\leqslant R_i+r_0+2\}
 \cap\bigcap_{i=1}^k\{\theta_{R_i}B\in G_n(x,r_0)\}.
\]
To pass to the limit in \(n\), we allow the left endpoints of the intervals to vary subject to closed constraints. For a compact target \(A\), let \(E_n(w)\) consist of the centers \(x\in A\) admitting intervals \(I_{i,v}=[t_{i,v},t_{i,v}+r_{|v|}]\), \(1\leqslant i\leqslant k\), \(v\in\{0,1\}^{\leqslant n}\). The initial intervals \(I_{i,\varnothing}\) lie in \([0,3k]\) in increasing order of \(i\), with gaps of at least \(1\) between successive intervals. For \(|v|<n\), the intervals \(I_{i,v0}\) and \(I_{i,v1}\) are contained in \(I_{i,v}\) in that order, with a gap of at least \(r_{|v|+1}\) between them. At every left endpoint we require \(|w(t_{i,v})-x|\leqslant r_{|v|}\).

\begin{lemma}\label{lem:fixed-forest}
Let \(A\subseteq\RR^2\) be nonempty and compact, \(z\in\RR^2\), \(a\geqslant a_0\), and \(k\in\NN\). There are \(\kappa,C>0\) such that the following hold for every \(n\geqslant0\).
\begin{enumerate}
\renewcommand{\labelenumi}{(\roman{enumi})}
\item The events \(F_{n,k}(x)\) are jointly Borel in the path and center, with
\[
 F_{n,k}(x)\in\FF_{3k},\quad F_{n+1,k}(x)\subseteq F_{n,k}(x).
\]
\item \(\PP_z(F_{n,k}(x))\geqslant\kappa\,p_n(r_0)^k>0,\quad x\in A.\)
\item \(\displaystyle \frac{\PP_z(F_{n,k}(x)\cap F_{n,k}(y))} {\PP_z(F_{n,k}(x))\PP_z(F_{n,k}(y))} \leqslant C|x-y|^{-\eta},\quad x,y\in A,\ x\ne y.\)
\item The sets \(E_n(w)\) are compact with closed graphs, and
\[
 E_{n+1}(w)\subseteq E_n(w),\quad
 \{x\in A:F_{n,k}(x)\text{ occurs}\}\subseteq E_n(w).
\]
\item \(\displaystyle\bigcap_{n\geqslant0}E_n(w) \subseteq\mathcal T_{k/(2a)}([0,3k]).\)
\end{enumerate}
\end{lemma}
\begin{proof}
We write \(p=p_n(r_0)\), \(T_i=R_i+r_0\), \(D_i=\{\theta_{R_i}B\in G_n(x,r_0)\}\), and
\[
 A_1=\{R_1\leqslant1\},\quad
 A_{i+1}=\{R_{i+1}\leqslant T_i+2\},\quad 1\leqslant i<k.
\]
Since \(D_i\in\FF_{T_i}\), \(A_i\in\FF_{R_i}\), and \(T_i\leqslant R_{i+1}\), we have \(F_{n,k}(x)\in\FF_{T_k}\). On this event,
\[
 T_i\leqslant1+(i-1)(r_0+2)+r_0<3i\leqslant3k,
 \quad 1\leqslant i\leqslant k.
\]
Thus \(F_{n,k}(x)\in\FF_{3k}\). Joint Borel measurability follows from that of the starting times and \(G_n\), while \(G_{n+1}\subseteq G_n\) gives the nesting. This proves (i).

For (ii), choose an open ball \(U\) with \(\dist(0,U)>r_0\). The Gaussian formula gives
\[
 \begin{aligned}
 v_0
 &=\inf_{x\in A}
   \PP_{z-x}\{B_{1/2}\in U,\ B_1\in\Ball(0,r_0/2)\}\\
 &=\frac1{\pi^2}\inf_{x\in A}
   \int_U\int_{\Ball(0,r_0/2)}
   e^{-|u-(z-x)|^2-|v-u|^2}\dd v\dd u\\
 &\geqslant
   \frac{r_0^2}{4\pi}\left(\int_U\dd u\right)
   \min_{x\in A,\ u\in\overline U,\ |v|\leqslant r_0/2}
   e^{-|u-(z-x)|^2-|v-u|^2}>0,
 \end{aligned}
\]
and
\[
 \begin{aligned}
 v_1
 &=\inf_{|y|\leqslant\sqrt{r_0}}
   \PP_y\{B_1\in U,\ B_2\in\Ball(0,r_0/2)\}\\
 &=\frac1{4\pi^2}\inf_{|y|\leqslant\sqrt{r_0}}
   \int_U\int_{\Ball(0,r_0/2)}
   e^{-(|u-y|^2+|v-u|^2)/2}\dd v\dd u\\
 &\geqslant
   \frac{r_0^2}{16\pi}\left(\int_U\dd u\right)
   \min_{|y|\leqslant\sqrt{r_0},\ u\in\overline U,\ |v|\leqslant r_0/2}
   e^{-(|u-y|^2+|v-u|^2)/2}>0.
 \end{aligned}
\]
After translation by \(x\), continuity forces a crossing of the entrance circle during \([1/2,1]\) in the first event. Thus \(R_1\leqslant1\) on this event, and \(\PP_z(A_1)\geqslant v_0\). For subsequent entrances, let
\[
 C_i=\{B_{T_i+1}-x\in U,\ B_{T_i+2}-x\in\Ball(0,r_0/2)\}.
\]
On \(C_i\), the path is outside the entrance circle at \(T_i+1\) and inside it at \(T_i+2\). Continuity therefore gives \(R_{i+1}\leqslant T_i+2\), so \(C_i\subseteq A_{i+1}\). By the strong Markov property at \(T_i\), the conditional law of \((B_{T_i+t}-x)_{t\geqslant0}\) given \(\FF_{T_i}\) is that of Brownian motion started at \(B_{T_i}-x\). On \(D_i\), the inclusion \(G_n(x,r_0)\subseteq G_0(x,r_0)\) gives \(|B_{T_i}-x|\leqslant\sqrt{r_0}\). Consequently,
\[
 \begin{aligned}
 \PP_z(A_{i+1}\mid\FF_{T_i})
 &\geqslant\PP_z(C_i\mid\FF_{T_i})\\
 &=\PP_{B_{T_i}-x}\{B_1\in U,\ B_2\in\Ball(0,r_0/2)\}\\
 &\geqslant\inf_{|y|\leqslant\sqrt{r_0}}
   \PP_y\{B_1\in U,\ B_2\in\Ball(0,r_0/2)\}=v_1.
 \end{aligned}
\]
At each entrance, \(|B_{R_i}-x|=r_0\), so \(\EE_z(\one_{D_i}\mid\FF_{R_i})=p\). Since \(A_{i+1}\in\FF_{R_{i+1}}\) and \(T_i\leqslant R_{i+1}\), conditioning first at \(R_{i+1}\) and then at \(T_i\) yields
\[
 \begin{aligned}
 \EE_z(\one_{A_{i+1}}\one_{D_{i+1}}\mid\FF_{T_i})
 &=\EE_z\!\left[\one_{A_{i+1}}
       \EE_z(\one_{D_{i+1}}\mid\FF_{R_{i+1}})\mid\FF_{T_i}\right]\\
 &=p\,\PP_z(A_{i+1}\mid\FF_{T_i})\geqslant v_1p,\quad\text{on }D_i.
 \end{aligned}
\]
For \(i=1\), \(F_{n,1}(x)=A_1\cap D_1\). Since \(A_1\in\FF_{R_1}\),
\[
 \begin{aligned}
 \PP_z(F_{n,1}(x))
 &=\EE_z\!\left[\one_{A_1}
      \EE_z(\one_{D_1}\mid\FF_{R_1})\right]\\
 &=p\,\PP_z(A_1)\geqslant v_0p.
 \end{aligned}
\]
For \(1\leqslant i<k\), we condition at \(T_i\) and use \(F_{n,i}(x)\in\FF_{T_i}\) and \(F_{n,i}(x)\subseteq D_i\) to obtain
\[
 \begin{aligned}
 \PP_z(F_{n,i+1}(x))
 &=\EE_z\!\left[\one_{F_{n,i}(x)}
      \EE_z(\one_{A_{i+1}}\one_{D_{i+1}}\mid\FF_{T_i})\right]\\
 &\geqslant v_1p\,\PP_z(F_{n,i}(x)).
 \end{aligned}
\]
Iterating gives
\[
 \PP_z(F_{n,k}(x))\geqslant v_0p\,(v_1p)^{k-1}
 =v_0v_1^{k-1}p_n(r_0)^k.
\]
This proves (ii) with \(\kappa=v_0v_1^{k-1}\).

To prove (iii), we compare the constructions at distinct centers \(x,y\in A\) and write
\[
 d=|x-y|,\quad \ell=\min\{j:2M\sqrt{r_j}<d\}.
\]
At a starting time \(S=S_{i,v}(x)\) with \(j=|v|\leqslant n\), we have \(|B_S-x|=r_j\). The strong Markov property gives
\[
 \PP_z\{\theta_SB\in G_{n-j}(x,r_j)\mid\FF_S\}
 =\PP_{B_S}\bigl(G_{n-j}(x,r_j)\bigr)=p_{n-j}(r_j).
\]
This event belongs to \(\FF_{S+r_j}\), as required by Lemma~\ref{lem:random-order-windows}. When \(n\geqslant\ell\), joint success supplies \(k2^\ell\) generation-\(\ell\) intervals at each center, with
\[
 B([S_{i,v}(x),S_{i,v}(x)+r_\ell])
 \subseteq\overline\Ball(x,M\sqrt{r_\ell}),\quad |v|=\ell,
\]
and the analogous inclusion at \(y\). The intervals at each center are separated, and the two confinement disks are disjoint by the choice of \(\ell\). Hence all \(2k2^\ell\) intervals have disjoint interiors. When \(n<\ell\), we retain only the \(k\) initial intervals at \(x\). Applying Lemma~\ref{lem:random-order-windows} in these two cases gives
\[
 \PP_z(F_{n,k}(x)\cap F_{n,k}(y))
 \leqslant
 \begin{cases}
 p_{n-\ell}(r_\ell)^{2k2^\ell},&n\geqslant\ell,\\
 p_n(r_0)^k,&n<\ell.
 \end{cases}
\]
Write \(\mathcal R_n(x,y)\) for the ratio in (iii). Dividing by the lower bound from (ii) and using \eqref{eq:iterate} and \eqref{eq:crude}, we obtain
\[
 \mathcal R_n(x,y)\leqslant\kappa^{-2}
 \begin{cases}
 q^{-2k(2^\ell-1)},&n\geqslant\ell,\\
 (qb)^{-k2^n},&n<\ell,
 \end{cases}
 \leqslant\kappa^{-2}K^{k2^\ell}.
\]
If \(\ell\geqslant1\), its minimality gives
\[
 d\leqslant2M\sqrt{r_{\ell-1}}=2M e^{-a2^{\ell-2}},\quad
 2^\ell\leqslant\frac4a\log\frac{2M}{d}.
\]
Since \(\eta=4k\log K/a\), it follows that
\[
 \mathcal R_n(x,y)
 \leqslant\kappa^{-2}K^{k2^\ell}
 \leqslant\kappa^{-2}(2M)^\eta d^{-\eta}.
\]
For \(\ell=0\), the estimate for \(n\geqslant\ell\) gives
\[
 \mathcal R_n(x,y)\leqslant\kappa^{-2}
 \leqslant\kappa^{-2}(1+\diam A)^\eta d^{-\eta}.
\]
We next prove (iv). Joint continuity of \((w,t)\mapsto w(t)\) makes the defining constraints closed in \(\mathcal C\times A\times[0,3k]^{k(2^{n+1}-1)}\). The time-parameter space is compact, so projection onto \((w,x)\) gives a closed graph with compact sections in \(A\). Restricting the family of intervals to \(|v|\leqslant n\) gives \(E_{n+1}(w)\subseteq E_n(w)\). On \(F_{n,k}(x)\), the left endpoints \(t_{i,v}=S_{i,v}(x)\) satisfy all the constraints, giving the other inclusion in (iv).

Finally, fix \(x\in\bigcap_nE_n(w)\). For each \(n\), choose admissible left endpoints \(t_{i,v}^{(n)}\) for \(|v|\leqslant n\). All endpoints lie in the compact interval \([0,3k]\). We extract convergent subsequences successively for the countably many labels \(i,v\); a diagonal choice then gives a single increasing sequence \(n_\ell\) such that, for every fixed \(i,v\),
\[
 t_{i,v}^{(n_\ell)}\longrightarrow t_{i,v},\qquad\ell\to\infty.
\]
Each constraint involves finitely many endpoints and is closed, so the limiting endpoints satisfy all the constraints. The resulting intervals define nonempty, pairwise separated compact sets
\[
 K_{x,i}=\bigcap_{j\geqslant0}\ \bigcup_{|v|=j}I_{i,v}.
\]
For \(t\in K_{x,i}\), choose \(|v_j|=j\) with \(t\in I_{i,v_j}\). Since \(|t-t_{i,v_j}|\leqslant r_j\), continuity gives
\[
 |w(t)-x|\leqslant|w(t)-w(t_{i,v_j})|+r_j\longrightarrow0.
\]
Thus \(K_{x,i}\subseteq L_x([0,3k])\).

We assign mass \(2^{-j}\) to each generation-\(j\) interval. The resulting measure \(\mu_i\) is the distribution of the unique limiting point under independent fair choices of successive children. Thus
\[
 \mu_i(I_{i,v})=2^{-|v|},\qquad
 \mu_i(K_{x,i})=\sum_{|v|=j}2^{-j}=1.
\]
For \(t\in K_{x,i}\), let \(I_{i,v_j}\) be the generation-\(j\) interval containing \(t\). Then
\[
 \mu_i(\{t\})\leqslant\mu_i(I_{i,v_j})=2^{-j}\longrightarrow0,
 \quad j\to\infty.
\]
Distinct generation-\(j\) intervals are separated by at least \(r_j\). Hence an interval \(J\) with \(0<|J|<r_0\) and \(r_j\leqslant|J|<r_{j-1}\) meets at most one interval of generation \(j-1\), and
\[
 \mu_i(J)\leqslant2^{-(j-1)}
 =2a\,h_0(r_j)\leqslant2a\,h_0(|J|).
\]
Every bounded subset of \(\RR\) lies in a closed interval of the same diameter, so we may compute \(\HH^{h_0}_\delta(K_{x,i})\) using interval covers. For \(0<\delta<r_0\) and any countable cover \(K_{x,i}\subseteq\bigcup_mJ_m\) with \(|J_m|<\delta\), the preceding estimates give
\[
 1=\mu_i(K_{x,i})\leqslant\sum_m\mu_i(J_m)
 \leqslant2a\sum_m h_0(|J_m|).
\]
Taking the infimum over these covers gives
\[
 \HH^{h_0}_\delta(K_{x,i})\geqslant\frac1{2a},\qquad
 \HH^{h_0}(K_{x,i})
 =\lim_{\delta\to0^+}\HH^{h_0}_\delta(K_{x,i})
 \geqslant\frac1{2a}.
\]
Since the \(K_{x,i}\) are pairwise separated,
\[
 \HH^{h_0}(L_x([0,3k]))
 \geqslant\HH^{h_0}\!\left(\bigcup_{i=1}^kK_{x,i}\right)
 =\sum_{i=1}^k\HH^{h_0}(K_{x,i})
 \geqslant\frac{k}{2a}.\qedhere
\]
\end{proof}

\subsection{The lower estimate for spatial targets}\label{sec:compactness}

Write \(C_1=2K_0\).

\begin{lemma}\label{lem:target-lower}
Let \(A\subseteq\RR^2\) be deterministic and compact with \(d=\dim_{\mathrm H}A>C_1m\) for some \(m>0\). For every \(T>0\), deterministic starting point \(z\), and \(0<s<d-C_1m\), there is positive probability that \(A\cap\mathcal T_m([0,T])\) contains a compact subset supporting a nonzero measure of finite \(s\)-energy.
\end{lemma}
\begin{proof}
We choose \(k\) so that
\[
 a=\frac{k}{2m}\geqslant a_0,\quad
 \eta=\frac{K_0k}{a}=C_1m.
\]
Then \(k/(2a)=m\), as required by Lemma~\ref{lem:fixed-forest}(v). We first prove the assertion for \(T=3k\). Since \(s+\eta<d\), Frostman's lemma gives a probability measure \(\sigma\) on \(A\) of finite \((s+\eta)\)-energy. In particular, \(\sigma\) is nonatomic, and
\[
 \mathcal E_\eta(\sigma)
 \leqslant(\diam A)^s\mathcal E_{s+\eta}(\sigma)<\infty.
\]
We use the successful centers to define the random measures
\[
 \nu_n(\dd x)=
 \frac{\one_{F_{n,k}(x)}}{\PP_z(F_{n,k}(x))}\,\sigma(\dd x).
\]
Lemma~\ref{lem:fixed-forest}(i),(ii) makes \(\nu_n\) a well-defined random finite measure. Tonelli's theorem gives
\[
 \EE_z\nu_n(A)
 =\int_A\frac{\PP_z(F_{n,k}(x))}{\PP_z(F_{n,k}(x))}\,\sigma(\dd x)=1.
\]
Since the diagonal has zero \(\sigma\otimes\sigma\)-measure, another application of Tonelli's theorem and Lemma~\ref{lem:fixed-forest}(iii) give
\begin{equation}\label{eq:target-lower-moments}
 \begin{aligned}
 \EE_z[\nu_n(A)^2]
 &=\iint_{A\times A}
   \frac{\PP_z(F_{n,k}(x)\cap F_{n,k}(y))}
        {\PP_z(F_{n,k}(x))\PP_z(F_{n,k}(y))}
   \,\sigma(\dd x)\sigma(\dd y)
 \leqslant C\mathcal E_\eta(\sigma),\\
 \EE_z\mathcal E_s(\nu_n)
 &=\iint_{A\times A}|x-y|^{-s}
   \frac{\PP_z(F_{n,k}(x)\cap F_{n,k}(y))}
        {\PP_z(F_{n,k}(x))\PP_z(F_{n,k}(y))}
   \,\sigma(\dd x)\sigma(\dd y)
 \leqslant C\mathcal E_{s+\eta}(\sigma).
 \end{aligned}
\end{equation}
We next obtain simultaneous bounds on the mass and energy with positive probability. Choose a constant \(C_2\geqslant1\) independent of \(n\) so that
\[
 \sup_n\EE_z[\nu_n(A)^2+\mathcal E_s(\nu_n)]\leqslant C_2.
\]
Paley--Zygmund gives
\[
 \PP_z\{\nu_n(A)\geqslant1/2\}
 \geqslant\frac{(\EE_z\nu_n(A))^2}{4\EE_z[\nu_n(A)^2]}
 \geqslant\frac1{4C_2}.
\]
With \(L=8C_2^2\), Markov's inequality gives
\[
 \PP_z\{\mathcal E_s(\nu_n)>L\}
 \leqslant\frac{\EE_z\mathcal E_s(\nu_n)}L
 \leqslant\frac1{8C_2}.
\]
Thus, for \(H_n=\{\nu_n(A)\geqslant1/2,\ \mathcal E_s(\nu_n)\leqslant L\}\),
\[
 \PP_z(H_n)\geqslant\frac1{4C_2}-\frac1{8C_2}
 =\frac1{8C_2}=:p>0,\quad n\geqslant0.
\]
Each union \(\bigcup_{n\geqslant N}H_n\) has probability at least \(p\), so continuity from above gives
\begin{equation}\label{eq:limit-event}
 \PP_z\left(\limsup_{n\to\infty}H_n\right)
 =\lim_{N\to\infty}\PP_z\left(\bigcup_{n\geqslant N}H_n\right)
 \geqslant p.
\end{equation}
On this event, we choose \(n_j\to\infty\) with \(w\in H_{n_j}\). For \(x,y\in A\), the energy kernel is at least \((\diam A)^{-s}\). Hence
\[
 (\diam A)^{-s}\nu_{n_j}(A)^2
 \leqslant\mathcal E_s(\nu_{n_j})\leqslant L.
\]
The measures have uniformly bounded total mass and are supported on the compact set \(A\). After passing to a subsequence, we may therefore assume that \(\nu_{n_j}\) converges weakly to a measure \(\nu\). Testing against the constant function \(1\) gives
\[
 \nu(A)=\lim_{j\to\infty}\nu_{n_j}(A)\geqslant1/2.
\]
For \(\ell\geqslant1\), define \(K_\ell(x,y)=\min\{\ell,|x-y|^{-s}\}\) for \(x\ne y\) and \(K_\ell(x,x)=\ell\). This kernel is bounded and continuous on \(A\times A\). Weak convergence of the product measures gives
\[
 \iint K_\ell(x,y)\,\nu(\dd x)\nu(\dd y)
 =\lim_{j\to\infty}
   \iint K_\ell(x,y)\,\nu_{n_j}(\dd x)\nu_{n_j}(\dd y)
 \leqslant L.
\]
Monotone convergence as \(\ell\to\infty\) yields \(\mathcal E_s(\nu)\leqslant L\).

We now verify that \(\nu\) has the required support. For each fixed \(N\), \(n_j\geqslant N\) for all sufficiently large \(j\). Lemma~\ref{lem:fixed-forest}(iv) and closedness of \(E_N(w)\) give
\[
 \operatorname{supp}\nu_{n_j}\subseteq E_{n_j}(w)\subseteq E_N(w).
\]
Since \(A\setminus E_N(w)\) is relatively open in \(A\), the Portmanteau theorem gives
\[
 \nu(A\setminus E_N(w))
 \leqslant\liminf_{j\to\infty}\nu_{n_j}(A\setminus E_N(w))=0.
\]
Thus \(\operatorname{supp}\nu\subseteq E_N(w)\) for every \(N\). Lemma~\ref{lem:fixed-forest}(v) now gives
\[
 \operatorname{supp}\nu\subseteq
 \bigcap_nE_n(w)\subseteq A\cap\mathcal T_m([0,3k]).
\]
The support of \(\nu\) is compact, proving the assertion for \(T=3k\).

We extend the result to arbitrary \(T>0\) by Brownian scaling. Set \(\lambda=\sqrt{3k/T}\) and \(\widetilde B_u=\lambda B_{u/\lambda^2}\). Then \(\widetilde B\) is Brownian motion started at \(\lambda z\). Writing \(\widetilde L_y\) for its level sets, we have
\[
 \widetilde L_{\lambda x}([0,3k])=\lambda^2L_x([0,T]),\quad
 \HH^{h_0}(\widetilde L_{\lambda x}([0,3k]))
 =\HH^{h_0}(L_x([0,T]))
\]
by Lemma~\ref{lem:gauge-properties}(iii). Apply the preceding construction to \(\widetilde B\) and \(\lambda A\), and push the resulting measure \(\widetilde\nu\) forward under \(y\mapsto y/\lambda\). The resulting nonzero measure \(\nu\) has compact support in \(A\cap\mathcal T_m([0,T])\) and satisfies
\[
 \mathcal E_s(\nu)=\lambda^s\mathcal E_s(\widetilde\nu)<\infty.\qedhere
\]
\end{proof}

\subsection{Uniform classification}\label{sec:classification}

\begin{proposition}\label{prop:maximum}
There exist a deterministic \(m_B\in(0,C_0]\) and an event of probability one on which
\[
 \sup_x\HH^{h_0}(L_x(I))=m_B
\]
for every \(I\in\mathcal I\).
\end{proposition}
\begin{proof}
We write \(M(I)=\sup_x\HH^{h_0}(L_x(I))\), which is measurable for bounded closed \(I\) by Lemma~\ref{lem:level-set-content}(ii). We first compare its laws on deterministic intervals. For \(0\leqslant a<b<\infty\), the process
\[
 W^{a,b}_u=\frac{B_{a+(b-a)u}-B_a}{\sqrt{b-a}},\quad 0\leqslant u\leqslant1,
\]
is standard Brownian motion. Since
\[
 L_{B_a+\sqrt{b-a}\,y}([a,b])
 =a+(b-a)\{u\in[0,1]:W^{a,b}_u=y\},
\]
Lemma~\ref{lem:gauge-properties}(iii) gives
\[
 M([a,b])
 =\sup_y\HH^{h_0}\{u\in[0,1]:W^{a,b}_u=y\}
 \overset{\mathrm d}=M([0,1]).
\]
Consequently, the random variables
\[
 X=M([0,1]),\quad X_1=M([0,1/2]),\quad X_2=M([1/2,1])
\]
have a common distribution function \(F\), and \(X_1,X_2\) are independent by independent increments. Since \(X\geqslant\max\{X_1,X_2\}\), we obtain
\[
 \begin{aligned}
 F(t)&=\PP(X\leqslant t)
 \leqslant\PP(X_1\leqslant t,\ X_2\leqslant t)\\
 &=\PP(X_1\leqslant t)\PP(X_2\leqslant t)
 =F(t)^2\leqslant F(t).
 \end{aligned}
\]
Thus \(F(t)\in\{0,1\}\) for every \(t\). Proposition~\ref{prop:upper} gives \(0\leqslant X\leqslant C_0\) almost surely. By right continuity of \(F\),
\[
 m_B:=\inf\{t:F(t)=1\}\in[0,C_0],\quad \PP(X=m_B)=1.
\]

To show that \(m_B>0\), we apply Lemma~\ref{lem:target-lower} to a closed disk \(A\), with \(T=1\) and \(0<m<2/C_1\). Since \(\dim_{\mathrm H}A=2>C_1m\), it gives
\[
 \PP(X\geqslant m)
 \geqslant\PP\{A\cap\mathcal T_m([0,1])\ne\varnothing\}>0.
\]
As \(X=m_B\) almost surely, this implies \(m_B\geqslant m>0\).

We now obtain a common event for all intervals. The equality in law above gives \(\PP\{M([a,b])=m_B\}=1\) for each deterministic \(0\leqslant a<b<\infty\). Hence
\[
 \Omega_M:=\bigcap_{a,b\in\mathbb Q,\,0\leqslant a<b}
       \{M([a,b])=m_B\},\quad \PP(\Omega_M)=1.
\]
On \(\Omega_M\), continuity from below for the increasing union \(L_x([0,\infty))=\bigcup_{N\geqslant1}L_x([0,N])\) gives
\[
 \begin{aligned}
 M([0,\infty))
 &=\sup_x\sup_{N\geqslant1}\HH^{h_0}(L_x([0,N]))\\
 &=\sup_{N\geqslant1}\sup_x\HH^{h_0}(L_x([0,N]))
 =\sup_{N\geqslant1}M([0,N])=m_B.
 \end{aligned}
\]
Every \(I\in\mathcal I\) contains a nondegenerate closed interval \(J\) with rational endpoints. Hence
\[
 m_B=M(J)\leqslant M(I)\leqslant M([0,\infty))=m_B.\qedhere
\]
\end{proof}

\begin{proof}[Proof of Theorem~\ref{thm:classification}]
We intersect the events of Propositions~\ref{prop:upper} and~\ref{prop:maximum}. Part~(ii) follows from Proposition~\ref{prop:maximum}. For~(iii), when \(0<\Lambda(h)<\infty\), Lemma~\ref{lem:gauge-properties}(ii) gives
\[
 m_B\Lambda(h)
 =\Lambda(h)\sup_x\HH^{h_0}(L_x(I))
 \leqslant\sup_x\HH^h(L_x(I))\leqslant C_0\Lambda(h).
\]
When \(\Lambda(h)=0\), Proposition~\ref{prop:upper} gives \(\HH^h(L_x(I))=0\) for every \(x\), which completes~(iii). For~(iv), \(m_B>0\) implies \(\mathcal S_0(I)\ne\varnothing\). If \(\Lambda(h)=\infty\), Lemma~\ref{lem:gauge-properties}(ii) gives \(\HH^h(L_x(I))=\infty\) for every \(x\in\mathcal S_0(I)\). Finally, (iii) and (iv) together prove~(i).
\end{proof}
\section{The spectrum of critical masses}\label{sec:mass-spectrum}

\subsection{An upper bound for spatial targets}

\begin{lemma}\label{lem:target-upper}
Let \(A\subseteq\RR^2\) be deterministic and analytic with \(d=\dim_{\mathrm H}A\), and let \(m>0\). Brownian motion from any deterministic starting point almost surely satisfies
\begin{equation}\label{eq:target-upper-dimension}
 \begin{aligned}
 \dim_{\mathrm H}\bigl(A\cap\mathcal T_m([0,\infty))\bigr)
 &\leqslant d-c_0m,\quad d\geqslant c_0m,\\
 A\cap\mathcal T_m([0,\infty))&=\varnothing,\quad d<c_0m.
 \end{aligned}
\end{equation}
\end{lemma}
\begin{proof}
We fix a starting point \(z\) and first suppose that \(A\) is bounded. Choose an integer \(R>|z|\). For a closed disk \(V\) of radius \(\varepsilon\), let \(N_V\) count the visits before \(\tau_R\) with delay \(\varepsilon\). Subsection~\ref{sec:upper} gives
\begin{equation}\label{eq:target-ball-tail}
 \PP_z\{N_V\geqslant u\log(1/\varepsilon)\}
 \leqslant C_{R,u}\varepsilon^{c_0u},\quad u>0
\end{equation}
for all sufficiently small \(\varepsilon\) uniformly over \(V\). Fix \(0<u<m\), and write
\[
 C_\rho(x)=\HH^{h_0}_\rho(L_x([0,\tau_R])),\quad
 A_N=\{x\in A:C_{2/N}(x)\geqslant u\},\quad N\geqslant6.
\]
Since \(C_{2/N}(x)\to\HH^{h_0}(L_x([0,\tau_R]))\) as \(N\to\infty\), we have
\[
 A\cap\mathcal T_m([0,\tau_R])\subseteq\bigcup_{N\geqslant6}A_N.
\]
When \(x\in A_N\cap V\) and \(\varepsilon<1/N\), the delayed-visit intervals and \(\{\tau_R\}\) cover \(L_x([0,\tau_R])\), so
\[
 u\leqslant C_{2/N}(x)\leqslant C_{2\varepsilon}(x)
 \leqslant N_Vh_0(\varepsilon).
\]
For \(s\geqslant0\) with \(s+c_0u>d\), choose deterministic covers \(\{V_{j,i}\}_i\) of \(A\) by closed disks of sufficiently small radii \(\varepsilon_{j,i}\) tending uniformly to zero, with \(\sum_i\varepsilon_{j,i}^{s+c_0u}\leqslant2^{-j}\). Their retained \(s\)-costs satisfy
\[
 Z_j=\sum_i(2\varepsilon_{j,i})^s
     \one_{\{N_{V_{j,i}}\geqslant u\log(1/\varepsilon_{j,i})\}},
 \quad \sum_j\EE_z Z_j\leqslant2^sC_{R,u}\sum_j2^{-j}<\infty.
\]
Tonelli's theorem gives \(\sum_jZ_j<\infty\) almost surely. On this event, each \(A_N\) is covered by the retained disks for all sufficiently large \(j\). For fixed \(N\), \(s>0\), and \(\delta>0\), these disks have diameter below \(\delta\) for all sufficiently large \(j\). Thus
\[
 \HH^s_\delta(A_N)\leqslant\liminf_{j\to\infty}Z_j=0,\quad
 \HH^s(A_N)=\lim_{\delta\to0^+}\HH^s_\delta(A_N)=0.
\]
When \(s=0\), the integer \(Z_j\) is eventually zero, so \(A_N=\varnothing\).

If \(d\geqslant c_0m\), each rational \(s>d-c_0m\) admits \(u\in(0,m)\) with \(s+c_0u>d\). Countable stability of Hausdorff dimension and intersection over these \(s\) give the stopped dimension bound. If \(d<c_0m\), we choose \(d/c_0<u<m\) and \(s=0\), giving \(A\cap\mathcal T_m([0,\tau_R])=\varnothing\).

To pass to the full path, continuity from below gives
\[
 A\cap\mathcal T_m([0,\infty))
 \subseteq\bigcup_{R\in\NN,\,R>|z|}
       \bigl(A\cap\mathcal T_v([0,\tau_R])\bigr),
 \quad0<v<m.
\]
We apply the stopped bounds to \(A\cap\overline\Ball(0,k)\) and intersect the events over integers \(k\geqslant1\), \(R>|z|\), and a fixed sequence \(v_j\in(0,m)\) with \(v_j\to m^-\). If \(d\geqslant c_0m\), the stopped bounds and countable stability give
\[
 \dim_{\mathrm H}\bigl(A\cap\mathcal T_m([0,\infty))\bigr)
 \leqslant d-c_0v_j,\quad j\geqslant1.
\]
Letting \(j\to\infty\) gives the dimension bound in \eqref{eq:target-upper-dimension}. If \(d<c_0m\), choose \(j\) with \(d<c_0v_j\). The stopped emptiness conclusion gives
\[
 A\cap\overline\Ball(0,k)\cap\mathcal T_{v_j}([0,\tau_R])
 =\varnothing,\quad k,R\in\NN,\ R>|z|.
\]
Taking the unions over \(k,R\) gives \(A\cap\mathcal T_m([0,\infty))=\varnothing\).
\end{proof}

\subsection{The deterministic spectrum}

\begin{proof}[Proof of Theorem~\ref{thm:mass-spectrum}]
For rational \(q\geqslant0\), Lemma~\ref{lem:level-set-content}(iv) and the half-interval argument of Proposition~\ref{prop:maximum} give a deterministic value \(d(q)\) for \(\dim_{\mathrm H}\mathcal S_q([0,1])\). Stationary increments and scaling give the same value on every nondegenerate bounded deterministic interval. We define
\[
 D(m)=\sup_{q\in\mathbb Q,\,q>m}d(q),\quad m\geqslant0.
\]
To obtain a common event, we use
\begin{equation}\label{eq:strict-threshold-continuity}
 \begin{aligned}
 \mathcal S_m([0,\infty))&=\bigcup_{N\geqslant1}\mathcal S_m([0,N]),\\
 \mathcal S_m(I)&=\bigcup_{q\in\mathbb Q,\,q>m}\mathcal S_q(I).
 \end{aligned}
\end{equation}
We intersect the event of Theorem~\ref{thm:classification} with the events above for rational thresholds and nondegenerate closed intervals with rational endpoints. On this event, \eqref{eq:strict-threshold-continuity} and any such interval \(J\subseteq I\) give
\[
 D(m)=\dim_{\mathrm H}\mathcal S_m(J)
 \leqslant\dim_{\mathrm H}\mathcal S_m(I)
 \leqslant\dim_{\mathrm H}\mathcal S_m([0,\infty))=D(m)
\]
for every \(m\geqslant0\) and \(I\in\mathcal I\). The definition makes \(D\) nonincreasing and right-continuous. Theorem~\ref{thm:classification} gives nonemptiness for \(m<m_B\) and emptiness for \(m\geqslant m_B\), so \(D(m_B)=0\). Since \(\mathcal T_m([0,1])\subseteq\mathcal S_0([0,1])\), Lemma~\ref{lem:target-lower} applied to a closed disk with \(T=1\) and the determinism of \(D(0)\) give
\[
 D(0)\geqslant2-C_1m,\quad0<m<2/C_1.
\]
Letting \(m\to0^+\) gives \(D(0)=2\).

The nonemptiness and dimension events for strict superlevel sets and their intersections are Borel by Lemma~\ref{lem:level-set-content}(iii),(iv) and bounded-time exhaustion. For the increment estimates, we condition on the past and apply the spatial-target bounds to the future path. Let \(\mathbb W\) be Wiener measure on \(\mathcal C_0=\{\gamma\in\mathcal C:\gamma(0)=0\}\). For \(t>0\), \(w\in C([0,t],\RR^2)\), and \(\gamma\in\mathcal C_0\), set
\[
 (w\otimes_t\gamma)(r)=
 \begin{cases}
  w(r),&0\leqslant r\leqslant t,\\
  w(t)+\gamma(r-t),&r>t.
 \end{cases}
\]
We first prove the lower increment bound. Fix \(0\leqslant m_1<m_2\leqslant m_B\), write \(v=m_2-m_1\), and take \(0<s<D(m_1)-C_1v\). For the Borel path event \(H_s=\{\dim_{\mathrm H}\mathcal S_{m_2}([0,2])\geqslant s\}\), define
\[
 p(w)=\int\one_{H_s}(w\otimes_1\gamma)\,\mathbb W(\dd\gamma),\quad
 \PP(H_s\mid\FF_1)=p(B|_{[0,1]})\quad\text{almost surely}.
\]
The function \(p\) is Borel. For almost every past path \(w\), Lemma~\ref{lem:level-set-content}(iii) and countable stability of Hausdorff dimension give a compact \(A\subseteq\mathcal S_{m_1}([0,1])\) with \(\dim_{\mathrm H}A>s+C_1v\). We apply Lemma~\ref{lem:target-lower} to this target and the future path started at \(w(1)\). With positive \(\mathbb W\)-probability, the path \(w\otimes_1\gamma\) then admits a compact \(K\) satisfying
\[
 K\subseteq A\cap\mathcal T_v([1,2])
 \subseteq\mathcal S_{m_2}([0,2]),\quad \dim_{\mathrm H}K\geqslant s.
\]
The second inclusion uses \(f_{[0,2]}(x)=f_{[0,1]}(x)+f_{[1,2]}(x)\), since singletons have zero \(h_0\)-measure. Thus \(p(w)>0\) for almost every past path \(w\). Since \(D(m_2)\) is deterministic, we obtain
\[
 \PP(H_s\mid\FF_1)>0\text{ almost surely},\quad
 \PP(H_s)>0,\quad D(m_2)\geqslant s.
\]
Let \(s\to(D(m_1)-C_1v)^-\) when \(D(m_1)>C_1v\), and use \(D(m_2)\geqslant0\) otherwise. Then
\begin{equation}\label{eq:spectrum-lower-increment}
 D(m_2)\geqslant D(m_1)-C_1v.
\end{equation}
By monotonicity, \(D\) is therefore \(C_1\)-Lipschitz on \([0,m_B]\).

For the upper increment bound, assume \(m_2<m_B\). We split the mass at a rational time. For \(x\in\mathcal S_{m_2}([0,\infty))\), write
\[
 \nu_x=\HH^{h_0}|_{L_x},\quad
 M_x=\nu_x([0,\infty))>m_1+v,\quad F_x(t)=\nu_x([0,t]).
\]
The measure \(\nu_x\) is finite and nonatomic. Hence \(F_x\) is continuous, \(F_x(0)=0\), and \(F_x(t)\to M_x\) as \(t\to\infty\). There is therefore a rational \(t>0\) such that
\[
 m_1<F_x(t)<M_x-v,\quad
 \nu_x([t,\infty))=M_x-F_x(t)>v.
\]
This gives
\begin{equation}\label{eq:mass-splitting}
 \mathcal S_{m_2}([0,\infty))
 \subseteq\bigcup_{t\in\mathbb Q_{>0}}
   \bigl(\mathcal S_{m_1}([0,t])\cap\mathcal S_v([t,\infty))\bigr).
\end{equation}
For rational \(t>0\), write
\[
 E_t=\mathcal S_{m_1}([0,t])\cap\mathcal S_v([t,\infty)),\quad
 d=D(m_1).
\]
The past target \(\mathcal S_{m_1}([0,t])\) is analytic and has dimension \(d\) almost surely. Since \(\mathcal S_v([t,\infty))\subseteq\mathcal T_v([t,\infty))\), Lemma~\ref{lem:target-upper} and the conditional probability formula above with \(\otimes_t\) give almost surely
\[
 \begin{cases}
  \PP(E_t\ne\varnothing\mid\FF_t)=0,&d<c_0v,\\[2pt]
  \PP(\dim_{\mathrm H}E_t>d-c_0v\mid\FF_t)=0,&d\geqslant c_0v.
 \end{cases}
\]
We intersect the resulting probability-one events over rational \(t>0\). The first case and \eqref{eq:mass-splitting} imply \(d\geqslant c_0v\), since \(\mathcal S_{m_2}([0,\infty))\ne\varnothing\). The second case and countable stability then yield
\begin{equation}\label{eq:spectrum-upper-increment}
 D(m_2)\leqslant D(m_1)-c_0v.
\end{equation}
Continuity of \(D\) extends \eqref{eq:spectrum-upper-increment} to \(m_2=m_B\).

The inclusions
\[
 \mathcal S_m(I)\subseteq\mathcal T_m(I)
 \subseteq\mathcal S_{m-\varepsilon}(I),\quad0<\varepsilon<m
\]
hold for \(0<m\leqslant m_B\). Continuity of \(D\) gives \(\dim_{\mathrm H}\mathcal T_m(I)=D(m)\). On the same event, an open \(U\) meeting \(B(I)\) contains \(B(J)\) for some nondegenerate interval \(J\subseteq I\). Then
\[
 \mathcal S_m(J)\subseteq\mathcal S_m(I)\cap U\subseteq\mathcal S_m(I),\quad
 \dim_{\mathrm H}\mathcal S_m(J)=\dim_{\mathrm H}\mathcal S_m(I)=D(m).\qedhere
\]
\end{proof}

\section{Exact masses}\label{sec:exact-mass}

We first compare the hitting events of exact-mass sets and mass superlevels, then use this comparison to recover their dimensions. Fix \(h\in\GG\) with \(\Lambda(h)<\infty\), \(m>0\), and a bounded deterministic interval \(J=[a,b]\) with \(a<b\). Write
\[
 f_J^h(x)=\HH^h(L_x(J)),\quad
 \mathcal X_m^h(J)=\{x:f_J^h(x)=m\},\quad
 \mathcal T_m^h(J)=\{x:f_J^h(x)\geqslant m\}.
\]
We record the accumulated mass at the current position by
\[
 Y_t=\HH^h\{r\in[a,t]:B_r=B_t\},\quad t\in J.
\]
Lemma~\ref{lem:level-set-content} makes the map \((t,w)\mapsto\HH^h\{r\in[a,t]:w(r)=w(t)\}\) Borel. Composing this map with the continuous adapted process \(t\mapsto(t,B_{\cdot\wedge t})\) shows that \(Y\) is optional. The hitting events below are analytic projections of Borel graphs; we work in the corresponding completed probability spaces.

\begin{lemma}\label{lem:exact-mass-hitting}
Let \(\eta\) be a random element of a standard Borel space \(\mathsf Z\), independent of \(B\). Suppose that \(R=R(\eta)\subseteq\RR^2\) and \(\{(z,x)\in\mathsf Z\times\RR^2:x\in R(z)\}\) is Borel. Up to a joint null set,
\[
 \{R\cap\mathcal X_m^h(J)\ne\varnothing\}
 =\{R\cap\mathcal T_m^h(J)\ne\varnothing\}.
\]
\end{lemma}
\begin{proof}
We use the usual augmentation of \((\FF_t\vee\sigma(\eta))_{t\geqslant0}\). Independence preserves the Brownian strong Markov property. Denote the two events by \(E\) and \(F\), and let \(H\) be the projection of
\[
 \mathcal O=\{(\omega,t):a\leqslant t\leqslant b,\ Y_t=m,\ B_t\in R\}.
\]
On the event of Proposition~\ref{prop:upper},
\begin{equation}\label{eq:exact-threshold-projection-inclusions}
 E\subseteq F\subseteq H.
\end{equation}
For \(x\in R\cap\mathcal T_m^h(J)\), write
\[
 t_x=\inf\{t\in J:\HH^h(L_x([a,t]))\geqslant m\}.
\]
Finiteness and nonatomicity of \(\HH^h|_{L_x(J)}\) give
\[
 t_x\in\operatorname{supp}(\HH^h|_{L_x(J)})\subseteq L_x(J),\quad
 Y_{t_x}=m.
\]

The set \(\mathcal O\) is optional. For \(\varepsilon>0\), the optional section theorem \cite[Theorem~IV.84]{DM78} gives a stopping time \(S\) with \(\PP(S<\infty)\geqslant\PP(H)-\varepsilon\) whose graph on \(\{S<\infty\}\) lies in \(\mathcal O\). On \(\{S<\infty\}\), the strong Markov property and point polarity imply almost surely that
\[
 B_t\ne B_S,\quad t>S;\quad f_J^h(B_S)=Y_S=m.
\]
Thus \(\{S<\infty\}\subseteq H\cap E\) modulo null sets, and
\[
 \PP(F\setminus E)\leqslant\PP(H\setminus E)
 \leqslant\PP(H)-\PP(S<\infty)\leqslant\varepsilon.\qedhere
\]
\end{proof}

To recover dimensions from Lemma~\ref{lem:exact-mass-hitting}, we use independent percolation test sets, following \cite[Lemma~3.4 and Remark~3.5]{KPX}.

\begin{lemma}\label{lem:hausdorff-probe}
For every \(s\in(0,2)\), there is a random \(F_\sigma\) set \(R_s\subseteq\RR^2\) with a Borel membership relation such that every deterministic analytic \(E\subseteq\RR^2\) satisfies
\[
 \begin{aligned}
 \PP(E\cap R_s\ne\varnothing)&=0,\quad \dim_{\mathrm H}E<s,\\
 \PP(E\cap R_s\ne\varnothing)&=1,\quad \dim_{\mathrm H}E>s.
 \end{aligned}
\]
\end{lemma}
\begin{proof}
The compact limit \(\Gamma_s\) of dyadic percolation in a fixed square with retention probability \(2^{-s}\) satisfies
\[
 \PP\{\Gamma_s\cap\Ball(x,r)\ne\varnothing\}\leqslant Cr^s.
\]
A ball meets a bounded number of comparable dyadic squares, each surviving with probability at most \(Cr^s\). For analytic \(E\) in this square, ball covers give
\[
 \PP(E\cap\Gamma_s\ne\varnothing)=0,\quad \dim_{\mathrm H}E<s.
\]
If \(\dim_{\mathrm H}E>s\), choose \(t\) with \(s<t<\dim_{\mathrm H}E\). The compact approximation theorem for analytic sets \cite{Howroyd} gives a compact \(K\subseteq E\) with
\[
 \HH^t(K)>0,\quad \dim_{\mathrm H}K\geqslant t>s.
\]
Hawkes' theorem \cite[Theorem~9.5(ii)]{MP} then yields
\[
 \PP(E\cap\Gamma_s\ne\varnothing)
 \geqslant\PP(K\cap\Gamma_s\ne\varnothing)>0.
\]
We take independent copies \(\Gamma_s^{k,j}\) in \([-k,k]^2\) and write \(R_s=\bigcup_{k,j\geqslant1}\Gamma_s^{k,j}\). If \(\dim_{\mathrm H}E>s\), choose \(k\) with \(p=\PP(E\cap\Gamma_s^{k,1}\ne\varnothing)>0\). Independence gives
\[
 \PP(E\cap R_s=\varnothing)\leqslant\lim_{N\to\infty}(1-p)^N=0.
\]
Each copy is an intersection of finite unions of retained closed squares, giving the Borel membership relation.
\end{proof}

To use these probes for random sets, let \(E(w)\subseteq\RR^2\) have a Borel membership relation, with \(w\) in a standard Borel space. For \(s\in\mathbb Q\cap(0,2)\), let \(\nu_s\) be the law of the parameter \(\eta\) defining \(R_s\), and set
\[
 H_s=\{(w,\eta):E(w)\cap R_s(\eta)\ne\varnothing\},\quad
 q_s(w)=\nu_s((H_s)_w).
\]
The set \(H_s\) is analytic. For any probability law \(\mu\) on the parameter space of \(w\), there are Borel sets \(H'_s,N_s\) such that
\[
 H_s\mathbin{\triangle}H'_s\subseteq N_s,\quad
 (\mu\otimes\nu_s)(N_s)=0.
\]
Applying Fubini's theorem to \(N_s\), we find that for \(\mu\)-almost every \(w\),
\[
 \nu_s((N_s)_w)=0,\quad
 q_s(w)=\int\one_{H'_s}(w,\eta)\,\nu_s(\dd\eta).
\]
Thus \(q_s\) is measurable in the completion of \(\mu\). Since \(\mu\) was arbitrary, \(q_s\) is universally measurable. We can now recover the dimension from Lemma~\ref{lem:hausdorff-probe}:
\[
 \dim_{\mathrm H}E(w)=
 \sup\bigl(\{s\in\mathbb Q\cap(0,2):q_s(w)>0\}\cup\{0\}\bigr).
\]
The dimension function is therefore universally measurable.

To preserve exact mass on intervals containing \(J\), we retain centers whose visits are confined to \(J\). Write
\[
 \mathcal C_m^h(J)=\mathcal X_m^h(J)\setminus
                  B\bigl([0,a)\cup(b,\infty)\bigr).
\]
Thus \(\HH^h(L_x(I))=m\) for every \(x\in\mathcal C_m^h(J)\) and every interval \(I\supseteq J\).

\begin{proposition}\label{thm:exact-mass-transfer}
For each fixed \(h\in\GG\) with \(\Lambda(h)<\infty\), \(m>0\), and bounded deterministic interval \(J=[a,b]\) with \(a<b\), there is an event of probability one on which
\begin{equation}\label{eq:exact-threshold-dimensions}
 \dim_{\mathrm H}\mathcal C_m^h(J)
 =\dim_{\mathrm H}\mathcal X_m^h(J)
 =\dim_{\mathrm H}\mathcal T_m^h(J).
\end{equation}
The events \(\{\mathcal X_m^h(J)\ne\varnothing\}\) and \(\{\mathcal T_m^h(J)\ne\varnothing\}\) agree up to a null set.
\end{proposition}
\begin{proof}
Taking \(R=\RR^2\) in Lemma~\ref{lem:exact-mass-hitting} proves the assertion about nonemptiness. With \(R_s\) independent of \(B\), the same lemma and Fubini's theorem on the completed product space give a single probability-one event on which
\[
 \PP_R(\mathcal X_m^h(J)\cap R_s\ne\varnothing)
 =\PP_R(\mathcal T_m^h(J)\cap R_s\ne\varnothing),\quad
 s\in\mathbb Q\cap(0,2).
\]
Lemma~\ref{lem:hausdorff-probe} therefore gives \(\dim_{\mathrm H}\mathcal X_m^h(J)=\dim_{\mathrm H}\mathcal T_m^h(J)\).

We next show that restricting to \(\mathcal C_m^h(J)\) preserves dimension. For this, we separate the path into the independent Brownian segments
\[
 \begin{aligned}
 U_r&=B_a-B_{a-r},&&0\leqslant r\leqslant a,\\
 W_t&=B_{a+t}-B_a,&&0\leqslant t\leqslant b-a,\\
 V_r&=B_{b+r}-B_b,&&r\geqslant0.
 \end{aligned}
\]
Write \(\mu_W\) for the law of \(W\) and \(\nu\) for the law of \((U,V)\). For fixed paths \(w,u,v\), set
\[
 \begin{aligned}
 E(w)&=\{x:\HH^h\{t\in[0,b-a]:w_t=x\}=m\},\\
 \mathcal R_w(u,v)&=-u((0,a])\cup
                 \bigl(w_{b-a}+v((0,\infty))\bigr).
 \end{aligned}
\]
With time measured from \(a\) and space from \(B_a\), we can write
\[
 \begin{aligned}
 E(W)&=\mathcal X_m^h(J)-B_a,\\
 \mathcal R_W(U,V)&=B\bigl([0,a)\cup(b,\infty)\bigr)-B_a,\\
 E(W)\setminus\mathcal R_W(U,V)&=\mathcal C_m^h(J)-B_a.
 \end{aligned}
\]
For fixed \(w\), point polarity applies to both parts of \(\mathcal R_w(u,v)\), so for every \(x\in\RR^2\),
\[
 \nu\{(u,v):x\in\mathcal R_w(u,v)\}=0.
\]
To obtain a dimension bound, let \(0<s<\dim_{\mathrm H}E(w)\) and use Frostman's lemma to choose a probability measure \(\sigma\) on a compact subset of \(E(w)\) with \(\mathcal E_s(\sigma)<\infty\). Integrating the preceding hitting probability against \(\sigma\), we obtain by Tonelli's theorem
\[
 \int\sigma(\mathcal R_w(u,v))\,\nu(\dd u,\dd v)
 =\int\nu\{(u,v):x\in\mathcal R_w(u,v)\}\,\sigma(\dd x)=0.
\]
Hence, for \(\nu\)-almost every \((u,v)\),
\[
 \begin{aligned}
 \sigma\bigl(E(w)\setminus\mathcal R_w(u,v)\bigr)&=1,\\
 \mathcal E_s\bigl(\sigma|_{E(w)\setminus\mathcal R_w(u,v)}\bigr)
 &=\mathcal E_s(\sigma)<\infty,\\
 \dim_{\mathrm H}\bigl(E(w)\setminus\mathcal R_w(u,v)\bigr)&\geqslant s.
 \end{aligned}
\]
To pass to the joint law, for \(s\in\mathbb Q\cap(0,2)\) consider the exceptional set
\[
 A_s=\{(w,u,v):\dim_{\mathrm H}E(w)>s,\quad
        \dim_{\mathrm H}(E(w)\setminus\mathcal R_w(u,v))<s\}.
\]
The graph of \(E\) is Borel by Lemma~\ref{lem:level-set-content}. The exterior traces are countable unions of compact path images, so \(\mathcal R_w(u,v)\) also has a Borel membership relation. The dimension measurability established above therefore makes \(A_s\) universally measurable. The bound for each fixed \(w\) and Fubini's theorem on the completed product space now imply
\[
 \nu((A_s)_w)=0\quad\text{for every }w,\qquad
 (\mu_W\otimes\nu)(A_s)
 =\int\nu((A_s)_w)\,\mu_W(\dd w)=0.
\]
Intersecting over rational \(s\) and translating back, we obtain
\[
 \dim_{\mathrm H}\mathcal C_m^h(J)
 =\dim_{\mathrm H}\mathcal X_m^h(J)
 \quad\text{almost surely},
\]
which proves \eqref{eq:exact-threshold-dimensions}.
\end{proof}

\begin{proof}[Proof of Theorem~\ref{thm:exact-mass-intro}]
We fix \(m\in(0,m_B)\) and intersect the event of Theorem~\ref{thm:mass-spectrum} with the events of Proposition~\ref{thm:exact-mass-transfer} for \(h=h_0\) and all nondegenerate closed intervals \(J\) with rational endpoints. For \(I\in\mathcal I\) and open \(U\) meeting \(B(I)\), we choose such a \(J\subseteq I\) with \(B(J)\subseteq U\).
\[
 \mathcal C_m^{h_0}(J)\subseteq\mathcal X_m(I)\cap U\subseteq\mathcal T_m(I),\quad
 \dim_{\mathrm H}\mathcal C_m^{h_0}(J)
 =\dim_{\mathrm H}\mathcal T_m(I)=D(m)>0.\qedhere
\]
\end{proof}

\end{document}